\documentclass[a4paper,twoside,11pt,reqno]{amsart}
\usepackage{amssymb}
\usepackage{amsmath}
\usepackage{amsfonts}

\newtheorem{theorem}{Theorem}

\title{    Hilbertian Interpolation}
\author{Marc Atteia} 
\email{marcatteia@orange.fr } 
\date{Version of \today}
\begin{document}
\maketitle 
\begin{abstract}  
I want to prove that all classical techniques of interpolation and approximation as Lagrange, Taylor,
Hermite interpolations Beziers interpolants, Quasi interpolants, Box splines and others (radial splines, simplicial splines) are derived from a \textbf{unique}
simple hilbertian scheme. For sake of simplicity, we shall consider only elementary examples which could be easily generalized.
\end{abstract}

 \section{ Hilbert spaces.}

  In the following, we say that a prehilbert space is a
vectorial space $\mathcal{H}$  with a scalar product $\left\langle
\cdot \mid \cdot \right\rangle $
 which is a bilinear and positive mapping from $\mathcal{%
H\times }$$\mathcal{H}$ into $%
\mathbb{R}
$ .
\\
  Then, we set $\left\Vert \cdot \right\Vert =\left(
\left\langle \cdot \mid \cdot \right\rangle \right) ^{\frac{1}{2}}$ ;$%
\left\Vert \cdot \right\Vert $ is a norm on $\mathcal{H}$ .
\\
 When $\mathcal{H}$  with the norm $\left\Vert \cdot \right\Vert $
is complete, one say that $\mathcal{H}$ is a Hilbert space.
\\
  So, we denote a (pre)hilbert space by the formula $\left(
\mathcal{H},\left\langle \cdot \mid \cdot \right\rangle \right) $ ( or $%
\left( \mathcal{H},\left\Vert \cdot \right\Vert \right) $).
\\
 \textbf{Examples :}
 \begin{align*}
H^{m}(a,b)&=\left\{ f\in C^{m-1}\left[ a,b\right] \text{ ; }\forall t\in %
\left[ a,b\right] \text{ , }f^{\left( m-1\right) }\left( t\right) = %
f^{\left( m-1\right) }\left( a\right) +\int_{a}^{t}\varphi \left( s\right)
ds\right .
\\
&
\left .\text{ with }\int_{a}^{b}\left\vert \varphi \left( s\right)
\right\vert ^{2}\text{ }ds\text{ }<+\infty \right\}
\end{align*}

 One can consider many equivalent topological structures on $%
H^{m}(a,b)$ relating to different scalar products,
 such that :%
\begin{eqnarray*}
\text{ If }f,g &\in &H^{m}(a,b)\text{ , } \\
(i)\text{ }\left\langle f\mid g\right\rangle &=&\text{}%
\int_{a}^{b}\sum_{j=0}^{m}f^{\left( j\right) }\left( t\right) .g^{\left(
j\right) }\left( t\right) \text{ }dt \\
(ii)\text{ }\left\langle f\mid g\right\rangle &=&\sum_{j=0}^{m-1}f^{\left(
j\right) }\left( a\right) .g^{\left( j\right) }\left( a\right) \text{}%
+\int_{a}^{b}f^{\left( m\right) }\left( t\right) .g^{\left( m\right) }\left(
t\right) \text{ }dt \\
(iii)\text{ }\left\langle f\mid g\right\rangle &=&\sum_{j=0}^{m-1}f\left(
\theta _{j}\right) .g\left( \theta _{j}\right) +\int_{a}^{b}f^{\left(
m\right) }\left( t\right) .g^{\left( m\right) }\left( t\right) \text{ }dt%
\text{ with }a<\theta _{0}<\theta _{1}<...<\theta _{m}<b
\end{eqnarray*}
 
\section{Reproducing hilbertian kernels}

 Let $\Omega $ an arbitrary set and $%
\mathbb{R}
^{\Omega }$ the set of mappings from $\Omega $ into $%
\mathbb{R}
$ .

 We say that :%
\begin{eqnarray*}
\mathcal{H} &\mathcal{\in }&\mathcal{H}ilb\mathcal{}\left(
\mathbb{R}
^{\Omega }\right) \text{ iff} \\
&&\mathcal{H}\text{ is a vectorial subspace of }%
\mathbb{R}
^{\Omega }\text{ and} \\
\forall t &\in &\Omega \text{ , }\exists C(t)>0\text{ such that }\forall
f\in \mathcal{H}\text{ , }\left\vert f\left( t\right) \right\vert \leq
C\left( t\right) \left\Vert f\right\Vert _{\mathcal{H}}
\end{eqnarray*}
\begin{theorem} 
 There exists an isomorphism between $\mathcal{H}ilb\mathcal{}%
\left(
\mathbb{R}
^{\Omega }\right) $ and $%
\mathbb{R}
_{+}^{\Omega \times \Omega }$.
 So, generally, it is equivalent to consider $\mathcal{H\in H}ilb%
\mathcal{}\left(
\mathbb{R}
^{\Omega }\right) $ 
 or its hilbertian kernel $H\in
\mathbb{R}
_{+}^{\Omega \times \Omega }$ except for numerical applications.
\end{theorem}
\section{ Examples of hilbertian kernels}

\subsection{ Polynomial kernels}

Let $m\in
\mathbb{N}
$ , $t_{0}\in $$%
\mathbb{R}
$ and :%
\[
\mathcal{P}_{m}\left( \
\mathbb{R}
\right) =\left\{
\begin{array}{c}
\text{polynomials }P\text{ on }\
\mathbb{R}
\text{ of degree, }d%
{{}^\circ}%
P\leq m\text{ } \\
\text{ with the scalar product } \\
\left\langle P\mid Q\right\rangle =\text{}\sum_{j=0}^{m}P^{\left( j\right)
}\left( t_{0}\right) .Q^{\left( j\right) }\left( t_{0}\right) \text{ }dt%
\end{array}%
\right\}
\]

Then :%
\[
\forall t,s\in
\mathbb{R}
\text{ , }H\left( t,s\right) =\sum_{j=0}^{m}\frac{\left( t-t_{0}\right)
^{j}}{j!}\frac{\left( s-t_{0}\right) ^{j}}{j!}
\]

\subsection{ A "Spline" kernel}

 $$\mathcal{H=}\left\{ f\in H^{1}\left( 0,1\right) \text{ };\text{ }%
f\left( 0\right) =0\text{ and }\forall f,g\in H^{1}\left( 0,1\right)
\text{ , }\left\langle f\mid g\right\rangle =\text{}\int_{0}^{1}f^{^{%
\prime }}\left( t\right) .g^{^{\prime }}\left( t\right) \text{ }dt\text{}%
\right\} $$

  Then : $H(t,s)=t-\left( t-s\right) _{+}$ .

\subsection{ Fourier's kernel}

  Let :%
\[
\mathcal{H=}\left\{
\begin{array}{c}
f\in H^{1}\left( 0,2\right) \text{ };\text{ }f\left( 0\right) =f\left(
2\right) \text{ , }\int_{0}^{1}f\left( t\right) \text{ }dt=0 \\
\text{and }\forall f,g\in H^{1}\left( 0,2\right) \text{ , }\left\langle
f\mid g\right\rangle =\text{}\int_{0}^{1}f^{^{\prime }}\left( t\right)
.g^{^{\prime }}\left( t\right) \text{ }dt\text{}%
\end{array}%
\right\}
\]

  Then :%
\[
H\left( t,s\right) =\sum_{k=1}^{\infty }\frac{\cos \left[ k\pi \left(
t-s\right) \right] }{\left( k\pi \right) ^{2}}
\]

  and :%
\[
\forall t\in \left[ \left( 0,2\right) \right] \text{ , }\forall f\in H%
\text{ , }f\left( t\right) =\left\langle f\mid H\left( \cdot ,t\right)
\right\rangle
\]

  That is the Fourier's development of $f$ which could be
extended periodically.
 
\section{Operations on hilbertian spaces / kernels}
 Let $\Omega _{1}$ and $\Omega _{2}$ be arbitrary sets and
$\mathcal{H}_{j}\in \mathcal{H}ilb\mathcal{}\left(
\mathbb{R}
^{\Omega _{j}}\right) $ $j=1,2$ .
\\
 We denote by $H_{j}\in
\mathbb{R}
_{+}^{\Omega _{j}\times \Omega _{j}}$ the hilbertian kernel of $%
\mathcal{H}_{j}$ , $j=1,2$ , and $\bot $ an operation
 such that $H=H_{1}\bot $ $H_{2}\in
\mathbb{R}
_{+}^{\Omega \times \Omega }$ , where $\Omega $ is a convenient set.
\\
 Then there exists $\mathcal{H\in }$ $%
\mathbb{R}
_{+}^{\Omega \times \Omega }$ , whose the hilbertian kernel is equal to$H.$
\\
 So, we set : $\mathcal{H=H}_{1}$ $\bot $ $\mathcal{H}_{2}$ .
\\
\textbf{Examples :}
\\
  (i) Let $H_{1},H_{2}\in
\mathbb{R}
_{+}^{\Omega \times \Omega }$ .

 Then :%
\[
\lambda H_{1}\text{ }\left( \lambda \geq 0\right) \text{ , }H_{1}\
+H_{2}\text{ , }\sup \left( H_{1},H_{2}\right) \text{ , }\inf \left( \
H_{1},H_{2}\right) \text{ , ...}\in
\mathbb{R}
_{+}^{\Omega \times \Omega }
\]

  (ii) Let $H_{j}\in
\mathbb{R}
_{+}^{\Omega _{j}\times \Omega _{j}}$ , $j=1,2$ .

 Then, if $\Omega =$$\Omega _{1}\times $$\Omega _{2}$ : \
\begin{eqnarray*}
H_{1}\otimes H_{2}\text{ , }\frac{1}{2!}\left( H_{1}\wedge H_{2}\right)
&=&\frac{1}{2!}\left( H_{1}\otimes H_{2}+H_{1}\otimes H_{2}\right) \text{, \
} \\
H_{1}\vee H_{2} &=&\frac{1}{2!}\left( H_{1}\otimes H_{2}-H_{1}\otimes
H_{2}\right) \in
\mathbb{R}
_{+}^{\Omega \times \Omega }\text{  }
\end{eqnarray*}
\section{  Extension}

  We can consider $Hilb\left( E\right) $ when $E$ is a
topological vectorial space and

 $ E$ $=\mathcal{C}^{0}\left( \Omega \right) $ or $E$ $=%
\mathcal{D}^{^{\prime }}\left( \Omega \right) $ $,....$.
 
\subsection{  Interpolation}

  Let us suppose that :

  $\mathcal{H}\in \mathcal{H}ilb\mathcal{}\left(
\mathbb{R}
^{\Omega }\right) $ , $\alpha _{0},$...$\alpha _{n}\in
\mathbb{R}
$ , and

$  \left( k_{0},\text{}k_{1},\text{.....,}k_{n}\right) $a
free system of vectors in $\mathcal{H}$ .

  We denote by $\Pi $ the following general interpolation
problem :%
\[
\Pi \text{ : }\inf \left\{ \left\Vert f\right\Vert ;\text{ }f\in \mathcal{H%
}\text{ and }\left\langle k_{j}\mid f\right\rangle \right\} =\alpha _{j}%
\text{ , }0\leq j\leq n
\]

  Then the solution of $\Pi $ is the element $\sigma
=Arg\left( \Pi \right) $ such that :%
\[
\forall t\in \Omega \text{ , }\sigma \left( t\right)
=\sum_{j=0}^{n}\lambda _{j}\text{}\left\langle k_{j}\mid H\left( \cdot
,t\right) \right\rangle \text{ with }\lambda _{j}\in
\mathbb{R}
\text{ , }0\leq j\leq n\text{}
\]
\\
  \textbf{Examples :}
\\
   \textit{(i)} \textit{Lagrange interpolation}.

  Let $\theta _{0},$...$\theta _{n}\in
\mathbb{R}
$ , $\mathcal{H=P}_{n}\left( \
\mathbb{R}
\right) $ , $\forall P,Q\in \mathcal{H}$ , $\left\langle P\mid
Q\right\rangle =\sum_{j=0}^{n}P\left( \theta _{j}\right) .Q\left( \theta
_{j}\right) $ .

  Then :
\[
H\left( s,t\right) =\sum_{j=0}^{n}L_{j,n}\left( s\right) .L_{j,n}\left(
t\right)
\]

  and :%
\[
\Pi \text{ : }\inf \left\{ \sum_{j=0}^{n}\left\vert f\left( \theta
_{j}\right) \right\vert ^{2};\text{ }f\in \mathcal{H}\text{ and }f\left(
\theta _{k}\right) =\alpha _{j}\text{ , }0\leq j\leq n\right\}
\]

  Thus :%
\[
\sigma =\sum_{j=0}^{n}\alpha _{j}L_{j,n}\left( =\sum_{j=0}^{n}\lambda
_{j}H\left( \cdot ,\theta _{j}\right) \right)
\]

   \textit{(ii) Taylor interpolation (Approximation).}

  Let : $t_{0}\in
\mathbb{R}
$ and :%
\[
\mathcal{H}=\left\{ f\text{ ; }\forall t\in
\mathbb{R}
\text{ , }f\left( t\right) =\sum_{j=0}^{\infty }\alpha _{j}\text{ }\frac{%
\left( t-t_{0}\right) ^{j}}{j!}\text{ with }\sum_{j=0}^{\infty }\text{ }%
\left\vert \alpha _{j}\right\vert ^{2}<+\infty \right\}
\]

\[
\text{(let us remark that}:\alpha _{j}=f^{\left( j\right) }\left( \
t_{0}\right) \text{)}
\]

  We suppose that :%
\[
\text{ }\forall f,g\in \mathcal{H}\text{ , }\left\langle f\mid
g\right\rangle =\sum_{j=0}^{\infty }f^{\left( j\right) }\left( \
t_{0}\right) \text{.}g^{\left( j\right) }\left( t_{0}\right)
\]

  Then $\mathcal{H}$ with $\left\langle \cdot \mid \cdot
\right\rangle $ is a Hilbert space with kernel $H$ such that :   \
\[
\forall s,t\in
\mathbb{R}
\text{ }H\left( s,t\right) =\sum_{j=0}^{\infty }\frac{\left(
s-t_{0}\right) ^{j}}{j!}\frac{\left( t-t_{0}\right) ^{j}}{j!}
\]

  So :%
\[
\Pi \text{ : }\inf \left\{ \sum_{j=0}^{\infty }\left\vert f^{\left(
j\right) }\left( t_{0}\right) \right\vert ^{2};\text{ }f\in \mathcal{H}%
\text{ and }f^{\left( k\right) }\left( t_{0}\right) =\alpha _{k}\text{
, }k\in
\mathbb{N}
\right\}
\]

  and :%
\[
\sigma =\sum_{j=0}^{\infty }\lambda _{j}\left( \frac{\partial ^{j}}{\partial
s^{j}}H\left( s,t\right) \right) _{s=t_{0}}=\sum_{j=0}^{\infty }\alpha _{j}%
\frac{\left( t-t_{0}\right) ^{j}}{j!}
\]

  (We note that : $\sigma ^{\left( j\right) }\left( \
t_{0}\right) =\alpha _{j}$) .

   \textit{(iii) Bezier- Bernstein interpolant.}

    $\left( \ast \right) $ Let us suppose that :
\begin{eqnarray*}
\mathcal{H}_{1} &=&\mathcal{H}_{2}=\mathcal{P}_{1}\left( \
\mathbb{R}
\right) \text{ , }\forall P,Q\in \mathcal{P}_{1}\left( \
\mathbb{R}
\right) \text{ , }\left\langle P\mid Q\right\rangle =P\left( 0\right)
.Q\left( 0\right) +P\left( 1\right) .Q\left( 1\right) \\
&&\text{ Then, the hilbrtian kernel of }\left( \mathcal{H}_{j}\text{ , }%
\left\langle \cdot \mid \cdot \right\rangle \right) \text{ is } \\
\forall s,t &\in &%
\mathbb{R}
\text{ ,}H_{j}\left( s,t\right) =st+\left( 1-s\right) \left( 1-t\right)
\text{ , }j=1,2\text{ .}
\end{eqnarray*}

  So, as in previous paragraphs, we have :%
\[
\Pi _{j}\text{ : }\inf \left\{ \left\Vert P\right\Vert _{j}^{2}\text{ };%
\text{ }P\in \mathcal{H}_{j}\text{ and }P\left( 0\right) =\alpha _{0}%
\text{ , }P\left( 1\right) =\alpha _{1}\right\}
\]

   and
\[
\forall t\in
\mathbb{R}
\text{ , }\sigma \left( t\right) =\alpha _{0}\left( 1-t\right) +\alpha _{1}t
\]

   \

$    \left( \ast \ast \right) $ \textit{Tensorial
product }

  Let $\mathcal{H}=\mathcal{H}_{1}\otimes \mathcal{H}_{2}=%
\mathcal{P}_{1}\left( \
\mathbb{R}
\times
\mathbb{R}
\right) $

  The hilbertian kernel $H$ of $\mathcal{H}$ with its
usual scalar product is such that :%
\[
\forall s,t,s^{\prime },t^{\prime }\in
\mathbb{R}
\text{ , }H\left( s,s^{\prime };t,t^{\prime }\right) =H_{1}\left(
s,t\right) .H_{2}\left( s^{\prime },t^{\prime }\right)
\]

  Moreover, we have :%
\begin{eqnarray*}
\Pi  &=&\Pi _{1}\otimes \Pi _{2}\text{ : } \\
&&\inf \left\{
\begin{array}{c}
\left\Vert P\right\Vert ^{2}\text{ };\text{ }P\in \mathcal{H}\text{ and }
\\
P\left( 0,0\right) =\alpha _{00}\text{ , }P\left( 1,0\right) =\alpha _{10}%
\text{ , }P\left( 0,1\right) =\alpha _{01}\text{ , }P\left( 1,1\right)
=\alpha _{11}\text{ }%
\end{array}%
\right\} \text{}
\end{eqnarray*}

   and :%
\[
\sigma \left( s,t\right) =\alpha _{00}\left( 1-s\right) \left( 1-t\right)
+\alpha _{01}\left( 1-s\right) t+\alpha _{10}\text{ }s\left( 1-t\right)
+\alpha _{11}\text{ }st\text{ }
\]

\bigskip

    $\left( \ast \ast \ast \right) $ \textit{%
Restriction to the diagonal}

  Now, we consider :%
\[
\forall s,t\in
\mathbb{R}
\text{ },\text{ }\widetilde{H}\left( s,t\right) =H_{1}\left( s,t\right)
.H_{2}\left( s,t\right)
\]

  Then, we can easily verify that $\widetilde{H}$ is a
hilbertian kernel.

  So, $\forall s\in
\mathbb{R}
$ , $\tau \left( s\right) =\sigma \left( s,s\right) $ is the solution of
the $\widetilde{\Pi }$ associated to $\widetilde{H}$ .

   \textit{One can, in the same way study box-splines
interpolants.}

   \textit{(iv) Polynomial spline interpolant. }

  Let $m\in
\mathbb{N}
$ , $m\geq 2$ and $\theta _{0},$...$\theta _{m-1}\in
\mathbb{R}
$ , $\theta _{0}<...<\theta _{m-1}$.

  Let :%
\begin{eqnarray*}
\mathcal{H} &=&\left\{ f\in H^{m}\left( a,b\right) \text{ ; }f\left( \theta
_{j}\right) =0\text{ , }j=0,1,...,\left( m-1\right) \right\}
\end{eqnarray*}
  with
the scalar product
$\left\langle \cdot \mid \cdot \right\rangle$ such that   
$$\forall f,g\in H^{m}\left( a,b\right) \text{ , }\left\langle f\mid
g\right\rangle =\int_{a}^{b}f^{\left( m\right) }\left( t\right) .g^{\left(
m\right) }\left( t\right) \text{ }dt
$$
If $H$ is the hilbertian kernel of $\mathcal{H}$ we
have $\forall s,t \in \left( a,b\right)$ :%
$$
H\left( s,t\right) =\left(
-1\right) ^{m}\left[
\begin{array}{c}
\displaystyle
G_{m}\left( s,t\right) -\sum_{j=0}^{m-1}L_{j,m-1}\left( s\right) G_{m}\left(
t,\theta _{j}\right) -\sum_{j=0}^{m-1}L_{k,m-1}\left( t\right) G_{m}\left(
\theta _{k},s\right) \\
\displaystyle
+\sum_{j,k=0}^{m-1}L_{j,m-1}\left( s\right) L_{k,m-1}\left( t\right)
G_{m}\left( \theta _{j},\theta _{k}\right)%
\end{array}%
\right]
$$ 
with 
$\displaystyle G_{m}\left( s,t\right) =\frac{1}{\left( 2m-1\right) !}\left(
s-t\right) _{+}^{2m-1}
$.
 \\
 Then we consider :
\[
\Pi \text{ : }\inf \left\{ \left\Vert f\right\Vert ;\text{ }f\in \mathcal{H%
}\text{ and }f\left( t_{j}\right) =\alpha _{j}\text{ , }0\leq j\leq n%
\text{ with }t_{0}<t_{1}<...<t_{n}\text{ , }n>m\right\}
\]
\\
  and :%
\[
\sigma =\sum_{j=0}^{n}\lambda _{j}H\left( \cdot ,t_{j}\right) \text{ is
called a polynomial spline of odd degree.}
\] 
\section{Approximation of the Dirac's functional}

 Let $\mathcal{H=H}^{2}\left( 0,1\right) $ and

  $\forall f,g\in \mathcal{H}$ , $\left\langle f\mid
g\right\rangle =f\left( 0\right) g\left( 0\right) +f^{\prime }\left(
0\right) g^{\prime }\left( 0\right) +\int_{0}^{1}f^{^{\prime \prime }}\left(
t\right) .g^{^{\prime \prime }}\left( t\right) $ $dt$ .

 Let :

\[
\mathcal{K=P}_{1}\left( 0,1\right) \text{ and }\mathcal{L=}\left\{
f\in H^{2}\left( 0,1\right) \text{ };\text{ }f\left( 0\right) =f^{\prime
}\left( 0\right) =0\right\}
\]

  and :

$\forall f,g\in \mathcal{K}$ , $\left\langle f\mid
g\right\rangle _{1}=f\left( 0\right) g\left( 0\right) +f^{\prime }\left(
0\right) g^{\prime }\left( 0\right) $ \

 $\forall f,g\in \mathcal{L}$ , $\left\langle f\mid
g\right\rangle _{2}=\int_{0}^{1}f^{^{\prime \prime }}\left( t\right)
.g^{^{\prime \prime }}\left( t\right) $ $dt$ .

 Then one can prove easily that :

 $\left( \mathcal{H},\left\langle \cdot \mid \cdot
\right\rangle \right) $is direct sum of the two Hilbert spaces $\mathcal{%
}\left( \mathcal{K},\left\langle \cdot \mid \cdot \right\rangle
_{1}\right) $ and \

$  \left( \mathcal{L},\left\langle \cdot \mid \cdot
\right\rangle _{2}\right) $ .

 Let $H$ , $K$ , $L$ the hilbertian kernels of $\left(
\mathcal{H},\left\langle \cdot \mid \cdot \right\rangle \right) ,$$%
\mathcal{}\left( \mathcal{K},\left\langle \cdot \mid \cdot \right\rangle
_{1}\right) ,$

 and $\left( \mathcal{L},\left\langle \cdot \mid \cdot \right\rangle
_{2}\right) $ (respectively).

  So,
\begin{eqnarray*}
H &=&K+L\text{ with: } \\
\forall s,t &\in &\left( 0,1\right) \text{ , }K\left( s,t\right) =1+st%
\text{ , }L\left( s,t\right) =\frac{\left( s-t\right) _{+}^{3}}{3!}+%
\frac{s^{2}t}{2!}-\frac{s^{3}}{3!}
\end{eqnarray*}

  Thus :%
\begin{eqnarray*}
\forall f &\in &\mathcal{H}\text{ , }f=f_{1}+f_{2}\text{ , and :} \\
\forall t &\in &\left( 0,1\right) \text{ , }f\left( t\right) =\left\langle
f_{1}\mid K\left( \cdot ,t\right) \right\rangle _{1}+\left\langle f_{2}\mid
L\left( \cdot ,t\right) \right\rangle _{2}\text{ } \\
&=&f\left( 0\right) +tf^{\prime }\left( 0\right) +\int_{0}^{1}\frac{\partial
^{2}L}{\partial s^{2}}\left( t,s\right) .f^{\prime \prime }\left( s\right)
\text{ }ds
\end{eqnarray*}

  We remark that : $\int_{0}^{1}\frac{\partial ^{2}L}{\partial
s^{2}}\left( t,s\right) .f^{\prime \prime }\left( s\right) $ $ds$ is the
Peano's kernel.

 Now, we set :%
\begin{eqnarray*}
E^{h}f &=&f\left( \cdot +h\right) \text{ , }D_{h}=\frac{1}{h}\left(
E^{h}-E^{0}\right) \text{ , }D_{h}^{2}=\text{}D_{h}\circ \text{}%
D_{h}\text{ and :} \\
L_{h}\left( s,t\right) &=&\text{}D_{h,1}^{2}\text{}D_{h,2}^{2}\text{ }%
L\left( s,t\right)
\end{eqnarray*}

  Then :

  $\left( \ast \right) $ $L_{h}$ is a B-spline and an
hilbertian kernel.

  $\left( \ast \ast \right) $ $\lim_{h\rightarrow
0}L_{h}\left( \cdot ,t\right) =\delta _{t}$ in $\left( \mathcal{C}%
^{0}\left( 0,1\right) \right) ^{\prime }$

    where $\delta _{t}$ is the Dirac's functional. 

   Let us remark that $\frac{\partial ^{2}L}{\partial s%
\text{ }\partial t}\left( t,s\right) =\delta _{t}\left( s\right) $ $\left(
=\delta _{s}\left( t\right) \right) $ in $\mathcal{D}^{\prime }\left(
0,1\right) $ only.

  So, we have :
\[
\forall t\in \left( 0,1\right) \text{ , }\forall f\in \mathcal{C}^{0}\left(
0,1\right) \text{ , }\lim_{h\rightarrow 0}\int_{0}^{1}L_{h}\left( s,t\right)
f\left( s\right) \text{ }ds=f\left( t\right) \text{ .}
\]
\\
\textbf{References. }
 \vspace{0.3cm}
   
    [1] Marc Atteia, Hilbertian kernels and spline functions, North Holland.
 \vspace{0.3cm}
 
    [2] Marc Atteia et Jean Gaches,
    Approximation hilbertienne, Splines, Ondelettes, Fractales,
    EDP Sciences - Grenoble Sciences.  
\end{document}